\documentclass[12pt]{article}
\usepackage{latexsym}                 
\usepackage{graphicx}
\usepackage{lmodern}
\usepackage[utf8,latin1]{inputenc}
\usepackage{latexsym}                 
\usepackage{color}                 
\usepackage{amssymb}
\usepackage{amsmath}
\usepackage{amsthm}
\usepackage{graphicx}
\usepackage{hyperref}
\usepackage{verbatim}
\usepackage{mathptmx}      
\usepackage{bm}
\usepackage{matlab-prettifier}
\usepackage[section]{placeins}
\def\dfrac#1#2{\displaystyle{\frac{#1}{#2}   }}
	
\def\fa{\hbox{ for all }}

\def\Xn{X_n=\{x_1,\ldots,x_n\}}

\newcommand{\R}{\ensuremath{\mathbb{R}}}

\def\bo1{\mathbf 1}
\def\bo0{\mathbf 0}

\def\calb{{\mathcal B}}

\def\calk{\mathcal{K}}

\def\c0{{\cal O}}

\def\biglf{\par\bigskip\noindent} 
\newcommand{\eq}{\end{equation}}
\def\M{Mat\'ern{ }\relax}
\def\eref#1{(\ref{#1}%
)}
\def\RSref#1{\ref{#1}%
}
\def\RSlabel#1{\label{#1}%
}
\def\RScite#1{\cite{#1}%
}

\newcommand{\bql}[1]{%
\begin{equation}\label{#1}%
}

\newtheorem{Definition}{Definition}
\newtheorem{Lemma}{Lemma}

\newtheorem{Theorem}{Theorem}
\newtheorem{Corollary}{Corollary}

\def\Wmdr{{\stackrel{\circ}{H}}{}_2^m(B_r(0))}
\def\Wmd{H_2^m(\R^d)}
\begin{document}
\begin{center}
  {\large \bf Optimal Compactly Supported Functions in Sobolev Spaces}
  \biglf
Robert Schaback\footnote{%
Institut für Numerische und Angewandte Mathematik,
Universit\"at G\"ottingen, Lotzestra\ss{}e 16--18, 37083 G\"ottingen, Germany,
{\tt schaback@math.uni-goettingen.de}}
\\
\biglf
Draft of \today
\end{center}
{\bf Abstract:}
This paper constructs unique compactly supported functions in Sobolev spaces
that have minimal norm, maximal support, and maximal central value, under
certain renormalizations. They may serve as optimized
basis functions in interpolation or
approximation, or as shape functions in meshless methods for PDE solving.
Their norm is useful for proving upper bounds for convergence rates
of interpolation in Sobolev spaces $H_2^m(\R^d)$,
and this paper gives the correct rate $m-d/2$ that arises
as convergence like $h^{m-d/2}$ for interpolation at meshwidth $h\to 0$
or a blow-up like $r^{-(m-d/2)}$ for norms of compactly supported
functions with support radius $r\to 0$. In Hilbert spaces
with infinitely smooth reproducing
kernels, like Gaussians or inverse multiquadrics,
there are no compactly supported functions at all, but in spaces with
limited smoothness, compactly supported functions exist and can be optimized
in the above way.
The construction is described in Hilbert space via projections,
and analytically via trace operators. Numerical examples are provided.\\
\biglf
MSC Classification:\\
65D12, 41A15, 46E35, 65D07, 47B32, 35A08, 35B07\\
Keywords:\\
Radial Basis Functions, Splines, bump functions, shape functions,
Sobolev Spaces
\section{Introduction}\RSlabel{SecIntro}
In general, a {\em test function} is a smooth compactly supported (CS)
function,
sometimes assumed to be of infinite smoothness,
and often called a {\em bump  function}. 
Such functions are very useful in Analysis at various places, e.g. as mollifiers
or as trial functions. 
\biglf
In reproducing kernel Hilbert spaces like global Sobolev spaces, they
also help to prove certain results, e.g. the optimal rates of
approximations of derivatives by scalable stencils
\RScite{davydov-schaback:2019-1}.
\biglf
But a closer look at  test functions in kernel-based spaces reveals that
they may not even exist in general. Consequently, this paper
focuses on test functions in kernel-based spaces, their existence,
their properties, and their applications. While writing, it emerged
that there seems to be no general theory of compactly supported functions
in Sobolev spaces, and this paper tries to fill the gap in that generality.
\biglf
First, the basics of kernel-based spaces are stated, for
fixing notation and the context. Section \RSref{SecBF} turns to
compactly supported functions, and Section \RSref{SecSimBuFuApp}
provides a simple application to the error analysis
of interpolation.
\biglf
Then Section \RSref{SecExPro} shows that compactly supported functions
may simply not exist at all, if Fourier transforms of kernels decay
exponentially. But in cases of algebraic decay, like for the \M{}
kernel generating Sobolev spaces, there exists unique norm-minimal
functions $b_r^*$ with support on the ball $B_r(0)$ of radius $r$
around the origin and value $b_r^*(0)=1$. By uniqueness, they are necessarily
radial.
\biglf
Other properties of these functions and their norms
are proven in Section \RSref{SecProp}. When norms are kept bounded, they have
the maximal possible value at zero and the least possible support radius. 
Section \RSref{SecBFSL} considers the behaviour of
the optimal norm $\|b_r^*\|$ as a function of $r$, proving 
$\|b_r^*\|=\Theta(r^{-m+d/2})$ for $r\to 0$
in case of Sobolev space $\Wmd$. This works by scaling arguments,
and it turns out that downscaling of optimal functions to a smaller
support radius loses optimality,
but still has asymptotically the same rate as above.
\biglf
Then Section \RSref{SecChar} gives a characterization in terms of Hilbert space
arguments via projections, and Section \RSref{SecAnalyChar} applies
standard Sobolev space trace arguments to get a computable
representation. Numerical examples are added in Section \RSref{SecExa},
while Section \RSref{SecConc} summarizes the results and raises
quite a number of issues that require further investigation.
\section{Basics of Spaces and Kernels}\RSlabel{SecBasics}
Let $K$ be a continuous symmetric positive definite
translation-invariant real-valued kernel on $\R^d\times \R^d$.
It generates a {\em native} Hilbert space $\calk$ of functions on $\R^d$
with inner product 
$(.,.)_K$ and the remarkable properties
$$
\begin{array}{rcll}
  K(x,y)
  &=&(K(x,\cdot),K(y,\cdot))_K &\fa x,\,y\in \R^d,\\
  f(x)&=&(f,K(x,\cdot))_K &\fa x\in \R^d,\,f\in \calk,\\
   K(x,y)
  &=&(\delta_x,\delta_y)_{\calk^*} &\fa x,\,y\in \R^d.\\ 
\end{array}
$$
All delta functionals $\delta_x\;:\;f\mapsto f(x)$ are continuous and have
a {\em kernel translate} $K(x,\cdot)$ as their Riesz representer.
The space $\calk$ is the closure of all kernel translates under the above inner
product. The theory of kernel-based spaces
is treated extensively in the books
\cite{buhmann:2003-1, wendland:2005-1,fasshauer-mccourt:2015-1}
and an earlier lecture note \RScite{schaback:1997-1}, but we recall 
the basic facts for convenience of readers and for fixing notation, as far
as needed in this paper.
\biglf
In Spatial Statistics, kernels arise as
covariance functions of a mean-zero random
field $R$ on $\R^d$ in the sense that each $x\in\R^d$ carries a
zero-mean second-order random variable $R(x)$ such that $K(x,y)=Cov(R(x),R(y))
\fa x,y\in\R^d$.
\biglf
In Real Analysis, Sobolev spaces $\Wmd$ for $m>d/2$ are Hilbert spaces
generated by {\em Whittle-\M} kernels
$$
K(x,y)=\|x-y\|_2^{m-d/2}K_{m-d/2}(\|x-y\|_2)
$$
up to factors,
where $K_\nu$ is the modified Bessel function of second kind. This case
is very important as well in Spatial Statistics, see
\RScite{porcu-et-al:2024-1} for an overview. We use the $\Wmd$ notation
instead of $W_2^m(\Omega)$, because we work on the full space and
define Sobolev spaces via Fourier transforms.
\biglf
The recovery of functions $f$ from values $f(x_1),\ldots,f(x_n)$ on a
set $\Xn\subset\R^d$ of data points can be uniquely done by a function
$s_{X_n,f}\in \calk$ generated by the kernel
translates $K(x_1,\cdot),\ldots,K(x_n,\cdot)$, and the
pointwise error bound is
\bql{eqErrBnd}
|f(x)-s_{X_n,f}(x)|\leq P_{X_n}(x)\|f\|_K \fa x\in \R^d,\,f\in\calk.
\eq
Here, the {\em Power Function} $P_X\;:\;\R^d\to \R$ arises. It is the
distance of $\delta_x$ to the span of all $\delta_{x_j},\;1\leq j\leq n$
in the norm of the dual $\calk^*$. In Spatial Statistics, $s_{X,f}(x)$
is the best unbiased linear predictor for $f(x)$ given all
$f(x_1),\ldots,f(x_n)$ under
the random field with covariance function $K$,
and $P_{X_n}^2(x)$ is the variance of this prediction. For what follows,
the Power Function has the alternative definition
\bql{eqPXn}
P_{X_n}(x)=\sup\{f(x)\;:\;f\in \calk,\;\|f\|_K\leq 1, \;f(x_j)=0,\;1\leq
j\leq n\}.
\eq
\section{Bump Functions}\RSlabel{SecBF}
We look at a special case of compactly supported functions first.
General compactly supported functions will come up later. 
\begin{Definition}\RSlabel{DefBF}
  A {\em bump function of support radius } $r$ is an element of the set
  $$
\displaystyle{\calb_r:=\{b\in \calk\;:\;b(0)=1, b(y)=0 \;\fa \;\|y\|_2\geq
  r\}   }
$$
of functions on $\calk$.
\end{Definition}
In spaces $\calk$ consisting of analytic functions, the set
$\calb$ will be empty, and Corollary \RSref{CorAnaly}
in Section \RSref{SecExPro} gives a sufficient criterion
in terms of kernel smoothness. But for most of the paper, existence of bump
functions is assumed, and then they exist for all scales, as shown in
Section \RSref{SecExPro}.
\biglf
This calls for {\em optimizing} bump functions, and a typical
application will be in Section \RSref{SecSimBuFuApp}.
\begin{Definition}\RSlabel{DefOBF}
A bump function $b_r^*\in \calb_r$ is {\em norm-minimal}, if
$$
\|b_r^*\|_K=\displaystyle{\inf \{\|b\|_{K }\;:\;b\in \calb_r\}
}=:\beta(r),  
$$
and we call $\beta\;:\;\R_{>0}\to \R_{>0}$ the {\em bump norm function}. 
\end{Definition}
These functions are the main topic of the paper, but
many results extend to general compactly supported functions. 
\section{A Simple Application}\RSlabel{SecSimBuFuApp}
The following
result is a good
motivation to investigate bump functions. Its consequences are
elaborated in \RScite{schaback:2023-1} concerning Trade-off Principles
between errors and stability.
Let $\Xn$ be a set of
data points in $\R^d$, and define the distance
  $$
dist(x,X_n)=\displaystyle{ \min_{1\leq j\leq n}
  \|x-x_j\|_2}
$$
from a point $x$ to $X_n$.
\begin{Theorem}\RSlabel{ThePovBnd}
Then
the Power Function $P_X$ satisfies
\bql{eqpowbnd}
P_X(x)\geq \dfrac{1}{\|b\|_{K }}
\eq
for all bump functions $b$ with support radius $dist(x,X_n)$
or larger. 
\end{Theorem} 
\begin{proof}
   Any such bump  function has zero
   as its interpolant, and thus
   $$
1=|b(x)|\leq P_X(x)\|b\|_{K }
$$
by the standard error bound \eref{eqErrBnd}.
\end{proof}
Inequalities like this provide lower bounds for
pointwise interpolation errors, and these bounds are best if
the norm is minimized. This is why we look at compactly supported functions with
minimal norm.
\biglf
We shall see in \eref{eqbetasob} that in Sobolev spaces $\Wmd$
 the bump norm function $\beta(r)$ behaves like
 $r^{d/2-m}$, and thus \eref{eqpowbnd} yields
 a simple counterpart to the standard upper bounds of the Power Function,
 proving their asymptotic optimality. An earlier but
 much more complicated proof goes back to \RScite{schaback:1995-1}.
 If the error is measured
 in $\calk$, all other interpolation techniques have larger errors.
 This proves that \eref{eqpowbnd} is also a lower bound for errors
 of all other interpolation processes in Sobolev spaces. The paper
 \RScite{davydov-schaback:2019-1} works in a somewhat different context,
 but it also uses bump functions to prove
 the optimal possible convergence rate for interpolation 
 and derivative approximation in Sobolev spaces.
 \biglf
 Papers using bump functions for different purposes are 
\RScite{narcowich-et-al:2006-1,demarchi-schaback:2010-1,%
davydov-schaback:2016-2, larsson-schaback:2023-1}, but there will be many
others. None of these papers look at bump functions
in detail. 
\biglf
 This calls for an investigation of  bump functions
 and the bump norm function for more general cases, but
 it turns out that this runs into serious unexpected difficulties
 that are interesting in their own right.
\section{Existence Problems}\RSlabel{SecExPro}
It is clear that compactly supported functions exist in all kernel--based spaces
that only require certain finite smoothness properties, like
Sobolev spaces $\Wmd$ for $m>d/2$.
Wendland functions are examples that even use the kernel itself
\RScite{wendland:1995-1,schaback:2010-1}.
\biglf
For other kernel-bases spaces, e.g. those based on multiquadrics and Gaussians,
the existence of compactly supported functions is a serious problem.
Clearly, there are $C^\infty$ functions with compact support,
but one has to find some that lie in the given native Hilbert space, i.e.
the Fourier transform has to satisfy a specific decay property.
The following negative result uses the fact that  
exponential decay of Fourier transforms implies local analyticity around zero,
contradicting compact support when applied to compactly supported functions.
\begin{Lemma}\RSlabel{LemExpDec}
Assume that the $d$-variate Fourier transform $\hat f$ of some
Fourier-transformable function $f$ on $\R^d$ has exponential decay
$$
|\hat f(\omega)|\leq C \exp(-c\|\omega\|_2) \fa \omega \in \R^d.
$$
Then the function is analytic in a ball around zero
of radius proportional to $c$, with a factor depending on $C$ and $d$.
\end{Lemma}
\begin{proof}
The derivatives of $f$ at zero are bounded by
$$
\begin{array}{rcl}
  |D^\alpha f(0)|
  &\leq &
  (2\pi)^{-d/2}\displaystyle{  \int_{\R^d}|\hat f(\omega)||(i\omega)^\alpha| d\omega}\\
  &\leq &
  (2\pi)^{-d/2}C\displaystyle{  \int_{\R^d}\exp(-c\|\omega\|_2)|\omega^\alpha|d\omega}\\
  &\leq &
  (2\pi)^{-d/2}C\displaystyle{  \int_{\R^d}\exp(-c\|\omega\|_1/\sqrt{d})|\omega|^{|\alpha|} d\omega}\\
&\leq &
  (2\pi)^{-d/2}C\displaystyle{
    \prod_{j=1}^d\int_\R\exp(-c|\omega_j|/\sqrt{d})|\omega_j|^{\alpha_j} d\omega_j}\\
&\leq &
  (2\pi)^{-d/2}C\displaystyle{2^d\left(\frac{\sqrt{d}}{c}\right)^{d+|\alpha|}\alpha!}\\
\end{array}
$$
due to
$$
\begin{array}{rcl}
  \displaystyle{ \int_\R\exp(-c|t|/\sqrt{d})|t|^{ n } dt}
  &=&
  2\displaystyle{ \int_0^\infty\exp(-ct/\sqrt{d})t^{ n } dt}\\
  &=&
  2\displaystyle{\frac{\sqrt{d}}{c} \int_0^\infty\exp(-s)\left(\frac{s\sqrt{d}}{c}\right)^{ n } ds}\\
  &=&
  2\displaystyle{\left(\frac{\sqrt{d}}{c}\right)^{ n +1}
    \int_0^\infty\exp(-s) s^ n  ds}\\
  &=&
  2\displaystyle{\left(\frac{\sqrt{d}}{c}\right)^{ n +1} n !}.
\end{array}
$$
This implies convergence of the Taylor series in a region
$$
\{x\in \R^d\;:\; |x_i|\leq r_i,\;1\leq i\leq d\}
$$
around zero 
defined by a vector $r=(r_1,\ldots,r_d)$ of {\em adjoint radii of
convergence}, if
$$
\overline{\lim}_{|\alpha|\to\infty}\left|\frac{D^\alpha f(0)}{\alpha!}r^\alpha\right|^{\frac{1}{|\alpha|}}\leq 1
  $$
holds \RScite{shabat:1976-1}. Up to $d$-dependent or constant factors, we have to look at
$$
\left|\frac{D^\alpha f(0)}{\alpha!}\right|^{1/|\alpha|}
\leq \left(\frac{\sqrt{d}}{c}\right)^{1+d/|\alpha|}
$$
to get the assertion.
\end{proof}
Gaussians admit arbitrarily large
$c$, while inverse multiquadrics have a fixed maximal $c$ depending on the
kernel parameters.
\begin{Theorem}\RSlabel{TheAnaly}
  The native space of the Gaussian consists
  of globally analytic functions, i.e. the local power series
  expansions all exist and converge globally.
  The native space of $d$-variate
inverse multiquadrics generated by the kernel $K(x,y)=(1+\|x-y\|_2^2)^{-m}$
consists of locally analytic functions, i.e. the local power series
expansions all exist and converge with at least a
fixed radius of convergence.\qed
\end{Theorem}
Applying analytic continuation if necessary, we get
\begin{Corollary}\RSlabel{CorAnaly}
  In spaces generated by kernels with at least exponentially decaying Fourier
  transform, there are no compactly supported functions. \qed
\end{Corollary}
Consequently, all  argumentations using compactly supported functions
fail for such cases.
\biglf
These results are not really
surprising. Recall that there are no compactly functions in univariate
complex analysis or in spaces of harmonic functions, by the Maximum Principle.
In addition, kernels arising from the Hausdorff-Bernstein-Widder
representation are analytic and therefore never compactly supported.
The above result is slightly different, but similarly negative.
\biglf
Other cases of Native Spaces without compactly supported
functions are those whose kernels
have power series  expansions in the domain of interest. These are handled
in \RScite{zwicknagl:2009-1,zwicknagl-schaback:2013-1}.
\biglf
We now prove existence of norm-minimal bump functions,
the generalization to compactly supported functions being evident later.
\begin{Theorem}\RSlabel{TheBufEx}
  If the set $\calb_r$ of bump functions of radius $r$ is not empty,
  the bump norm function $\beta(r)$ is attained at  a unique bump function
  $b_r^*\in \calb_r$, i.e. $\beta(r)=\|b^*_r\|_{K }$.
\end{Theorem} 
\begin{proof} This follows from a standard variational argument.
  We start with an admissible bump function $\tilde b_r$ and approximate it from
  the closed linear subspace  
  $$
\displaystyle{\calb^0_r:=\{b\in \calk\;:\;b(0)=0,\;b(y)=0 \;\fa \;\|y\|_2\geq
  r\}   }
$$
to get some $\tilde b_r^0\in \calb^0_r$ with the orthogonality property
$$
(\tilde b_r-\tilde b_r^0,b_r^0)_K=0 \fa b_r^0 \in \calb^0_r.
$$
Then we define $b_r^*:=\tilde b_r-\tilde b_r^0\in \calb_r$ and take an
arbitrary $b_r\in \calb_r$ to get
$$
\begin{array}{rcl}
  \|b_r\|_K^2
  &=&
  \|b_r-b_r^*+b_r^*\|_K^2\\
  &=&
  \|b_r-b_r^*\|_K^2+2(b_r-b_r^*,b_r^*)_K +\|b_r^*\|_K^2\\
  &=&
  \|b_r-b_r^*\|_K^2+\|b_r^*\|_K^2\\
  &\geq & \|b_r^*\|_K^2.\\
\end{array}
$$
\end{proof}
The proof generalizes to any compact domain $\Omega\subset\R^d$
and an arbitrary
point $x$ in its interior for fixing the value 1.
Uniqueness follows similarly, and implies
\begin{Corollary}\RSlabel{CorBufEx}
  For {\em Radial Basis Functions} (RBFs), norm-minimal
  centralized bump functions
  on balls are radial, if they exist.
\qed
\end{Corollary}
This makes it easy to deal with such functions, once the radial
form is known or calcutated to sufficient precision. 
\section{Properties of Bump Functions}\RSlabel{SecProp}
From the definition, we have
\begin{Theorem}\RSlabel{Themono}
  Norms $\|b_r^*\|_K=\beta(r)$ of optimal bump functions $b_r^*$ increase
  with decreasing radius. 
\end{Theorem}
A second optimality principle is
\begin{Theorem}\RSlabel{Theoptrad}
Norm-minimal bump functions have maximal support radius under all
bump functions with norm up to one. 
\end{Theorem} 
\begin{proof}
  Assume a bump function $b_\rho$ for radius $\rho<r$.
  Then
  $$
  \|b_\rho\|_K \geq \beta(\rho)\geq \beta(r)=\|b_r^*\|_K.
  $$
\end{proof} 
The bump norm function cannot decrease to zero for large $r$,
because there is a positive lower bound:
\begin{Lemma}\RSlabel{Lemlow}
  $$
\beta^2(r)\geq 1/K(0,0) \fa r>0.
  $$
\end{Lemma}
\begin{proof}
  Let $b$ be any bump function at $0$, of arbitrary radius $r$. Then
  $$
1=b(0)=(b,K(0,\cdot))_{K }\leq \|b\|_{K }\|K(0,\cdot)\|_{K }=\|b\|_{K }\sqrt{K(0,0)}.
  $$
\end{proof}
This is just a special case of the general embedding inequality
$$
|f(x)|=|(f,K(x-\cdot)|\leq \|f\|_K \|K(x-\cdot)\|_K=\|f\|_K
\sqrt{K(0,0)} \fa f\in \calk,\;x\in \R^d.
$$
A third optimality is
\begin{Theorem}\RSlabel{TheValopt}
Under all bump functions with norm up to one, the maximum value at zero is
attained for $b_r^*/\|b_r^*\|_K$. 
\end{Theorem}
\begin{proof}
The optimal value is surely not less than $1/\|b_r^*\|_K$. If $b_r$ is any bump
function
with $\|b_r\|_K\leq 1$, the function $b_r/b_r(0)$ is admissible for norm
minimality, and thus
$$
\dfrac{1}{b_r(0)}\geq \dfrac{\|b_r\|_K}{b_r(0)}\geq \|b_r^*\|_K
$$
proves that the optimal value is at most  $1/\|b_r^*\|_K$.
\end{proof}
We can rewrite this as
$$
\dfrac{1}{\|b_r^*\|_K}=\sup\{f(0)\;:\;\|f\|_K\leq 1,\;f(x)=0 \fa  \|x\|_2\geq r\}
$$
and compare with \eref{eqPXn} to get the Power Function value 
$P_{\R^d\setminus B_r(0)}(0)$ for transfinite interpolation on all points
outside the ball $B_r(0)$. Note that this realizes the optimal
case in \eref{eqpowbnd}.
\section{Scaling Laws}\RSlabel{SecBFSL}
Now we check bump functions under scaling. To avoid clashes of indices,
we use a scaling operator $S_r$ acting on functions $f$ on $\R^d$
as $S_r(f)(x)=f(x/r) \fa x\in \R^d$ and $r>0$. 
If $f$ is supported on a ball $B_1(0)$ with radius $r$ around the origin,
the scaled function $S_r(f)$
is supported on the ball $B_{r}(0)$.
\biglf
But norm-minimality does not scale that way. If the kernel $K$ is fixed,
the optimal bump function $b^*_1$ for radius $1$ may be scaled
into $S_1 b^*_1$ to be admissible on $B_r(0)$,
but this need not be equal to $b^*_r$.
By norm-minimality,
$$
\|b^*_{r}\|_K\leq \|S_r b^*_1\|_K,
$$
and this is all we know. Figure \RSref{Fig1Dbumps} below demonstrates
numerically that these two differ, because the shape
of $b_r^*$ changes with $r$ in a nontrivial way.
\biglf
But scaling functions is also related to scaling kernels.
Therefore we locally change the notation $b_r^*$ to $b_{r,K}^*$ if the kernel $K$ is
used, and we may scale kernels $K$ with $S_r$ into $S_rK$ as well. 
However, since $\|f\|_{K} = \|S_{r}f\|_{S_{r}K}$ holds
by \RScite{larsson-schaback:2023-1} for all
functions in the native space of $K$, we
can apply this to all non-optimal bump functions of the appropriate
support radii to get 
$$
\|S_{r}b^*_{1,K}\|_{S_{r}K}=\|b^*_{1,K}\|_K=\inf \|b_{1}\|_K =\inf \|S_r b_{1}\|_{S_r K}=\inf
\|b_{r}\|_{S_r K}=\|b^*_{r,S_rK}\|_{S_r K}.
$$
Therefore the law 
$$
S_{r}b^*_{1,K}=b^*_{r,S_rK} 
$$
connects scaling of norm-minimal bump functions with scaling of the kernel.
\biglf
But it is more interesting to fix the kernel and vary the radius $r$.
We focus on kernels $K$ with
\bql{eqKomega}
\hat K(\omega)=\Theta(1+\|\omega\|_2^{-\beta}),\;\omega\in \R^d
\eq
with $\beta >d$ and a positive constant $C$. For the \M kernel
generating Sobolev space
$\Wmd$ we have $\beta=2m>d$.
\biglf
We bound the $d$--variate
Fourier transform of any scaled function $S_rf(x)=f(x/r)$
via $\widehat{S_rf}(\omega)=r^d\hat f(r\omega)$ to get 
$$
\begin{array}{rcl}
  \|S_rf\|_K^2
  &=&
  \displaystyle{\int_{\R^d} |\widehat{S_rf}(\omega)|^2(1+\|\omega\|_2^\beta) d\omega   }\\ 
  &=&
  \displaystyle{r^{2d}\int_{\R^d} |\hat f(r\omega)|^2(1+r^{-\beta}\|r\omega\|_2^\beta) d\omega   }\\ 
  &=&
  \displaystyle{r^{d-\beta}\int_{\R^d} |\hat f(\eta)|^2(r^\beta+\|\eta\|_2^\beta) d\eta   }. 
\end{array}
$$
up to constant factors. Now for $r\leq 1$ we find
$$
\begin{array}{rcl}
  \|S_rf\|_K^2
  &\leq & r^{d-\beta}\|f\|_K^2,
\end{array}
$$
and for lower bounds we use
$$
\begin{array}{rcl}
  \|S_rf\|_K^2
  &\geq &
  r^{d-\beta}|f|_K^2\\
\end{array}
$$
using the correspondent 
seminorm.  
When we apply all of this to bump functions, the seminorm
can be bounded below by the norm up to a factor, due to Poincar\'e inequalities.
\begin{Theorem}\RSlabel{TheSobBNF}
  In spaces generated by kernels with \eref{eqKomega},
  the bump norm function behaves like
\bql{eqbetasob}
  \beta(r) =\|b_r^*\|_K=\Theta(r^{d/2-\beta/2}) \hbox{ for } r\to 0.
  \eq
  In Sobolev space $\Wmd$, this holds for $\beta=2m$. 
\end{Theorem}
\begin{proof}
  For an upper bound we use
  $$
\|b_r^*\|^2_K\leq \inf \|b_r\|^2_K= \inf \|S_r b_1\|^2_K\leq r^{d-\beta} \|b_1\|^2_K 
$$
inserting arbitrary bump functions $b_1$ and $b_r$, The other direction is
$$
\begin{array}{rcl}
\|b_r^*\|^2_K=\|S_r S_{1/r} b_r^*\|^2_K&\geq&  r^{d-\beta}|S_{1/r} b_r^*|^2_K\\
&\geq&  r^{d-\beta}|S_{1/r} b_r^*|^2_K \;\inf |b_1|_K^2\\
&\geq& C\;\inf \|b_1\|_K^2=C\|b_1^*\|^2_K
\end{array} 
$$
\end{proof} 
We now  check scaled versions of norm-minimal bump functions. 
Our main tool is 
\begin{Theorem}\RSlabel{TheELRS}\cite[Thm. 1]{larsson-schaback:2023-1}\\
For all $f$ in the native space of kernels with \eref{eqKomega} and all $\epsilon >0$,
$$
\|S_{1/\epsilon} f\|^2_{K}=\|f\|^2_{S_\epsilon K}=
\Theta\left(\epsilon^d \max(1,\epsilon^{-\beta})   \right)\|f\|^2_{K}.\qed 
$$
\end{Theorem} 
We apply this to norm-minimal bump functions
and get
$$
\|b_{r}^*\|^2_{K}\leq \|S_{1/r}b^*_{1}\|^2_{K}=
\Theta\left(r^d \max(1,r^{-\beta})   \right)\|b^*_1\|^2_{K}.
$$
Therefore we can scale norm-minimal bump functions without asymptotic loss:
\begin{Corollary}\RSlabel{CorELRSscale}
  Under the assumption \eref{eqKomega}, scaled optimal bump
  functions satisfy
$$
\|S_{1/\epsilon} b_r^*\|^2_{K}=
\Theta\left(\epsilon^d \max(1,\epsilon^{-\beta})   \right)\|b_r^*\|^2_{K}, 
$$
where $S_{1/\epsilon} b_r^*$ is a
nonoptimal bump function with support radius $r\epsilon$.
\end{Corollary}
Figure \RSref{Fig1Dbumps} shows that the true bump functions
$b_r^*$ differ from scaled versions $b_1^*(\cdot/r)$, because their shape varies
nontrivially with $r$. 
\biglf
Again, all of this will generalize to compactly supported functions
on arbitrary domains $\Omega$ with a fixed interior point $x$.
By a shift, $x$ can be assumed to be the origin, and then the domain
is scaled as $S_r\Omega$.
\section{Characterization in Hilbert Space}\RSlabel{SecChar}
To get a constructive characterization of bump functions,
consider the closed subspace
$$
V_r:=\{ v \in \calk\;:\; v(x)=0 \fa  \|x\|\geq r\}
$$
of $\calk$, and there are bump functions iff
\bql{eqKvK0}
v(0)=0 \fa v\in V_r
\eq
is not satisfied.
\begin{Theorem}\RSlabel{TheChar}
  If there are bump functions at all, the unique norm-minimal  bump function
  $b_r^*$ on $B_r(0)$ has the form
\bql{eqbrgr}
b_r^*:=g_r/g_r(0)
\eq
for the projection $g_r$ of $K(0,\cdot)$ onto $V_r$.
The squared norm of the solution is 
\bql{eqrnnorm}
\|b_r^*\|_K^2=\dfrac{\|g_r\|_K^2}{g_r^2(0)}=\dfrac{1}{g_r(0)}.
\eq
\end{Theorem} 
\begin{proof}
  If $K(0,\cdot)$ were orthogonal to $V_r$, the equations
  $$
(K(0,\cdot),v_r)_K=0=v_r(0)
  $$
  would hold for all $v_r\in V_r$, implying \eref{eqKvK0}
  and the nonexistence of bump functions.
Let $P_{r}$ be the Hilbert space projector onto $V_{r}$. Then
$g_r:=P_{r}K(0,\cdot)$ is uniquely defined and nonzero.
Furthermore,
$K(0,\cdot)-g_r$ is orthogonal to $V_r$, i.e.
\bql{eqvrgrvr}
\begin{array}{rcl}
(K(0,\cdot)-g_r,v_r)_K &=&0,\\
v_r(0)&=& (g_r,v_r)_K 
\end{array} 
\eq
for all $v_r\in V_r$.
In particular, setting $v_r=g_r$ yields 
$$
g_r(0)= (g_r,g_r)_K>0.
$$
Therefore \eref{eqbrgr}
solves the problem and \eref{eqrnnorm} holds.
In fact, for any other bump function $b_r$
the difference $w_r=b_r^*-b_r$ is in $V_r$ with $w_r(0)=0$,
and
$$
(b_r^*,w_r)=\dfrac{1}{g_r(0)}(g_r,w_r)_K=\dfrac{w_r(0)}{g_r(0)}=0
$$
proves
$$
\|b_r\|^2_K= \|b_r^*-w_r\|^2_K=\|b_r^*\|^2_K+\|w_r\|^2_K
\geq \|b_r^*\|^2_K.
$$
\end{proof} 
By the second identity of \eref{eqvrgrvr}, $g_r$ is the Riesz representer of
the functional $\delta_0$ on the Hilbert space $V_r$.
We can generalize this to $g_{r,x}:=P_{r}K(x,\cdot)$ to get the Riesz
representer of $\delta_x$, and the kernel
\bql{eqKr}
K_r(x,y):=(g_{r,x},g_{r,y})_K \fa x,y\in B_r(0)
\eq
is reproducing on the Hilbert space $V_r$. This kernel lives
on the interior of the ball $B_r(0)$ centred at zero with radius $r$.
Figure \RSref{Fig1Dgrx} shows the kernel translates $g_{1,x}$
in Sobolev space $\stackrel{\circ}{W}{}_2^1(-1,1))$. The calculations
are based on
Section \RSref{SecChar}.
\biglf
\def\RSh{6.0cm}
\def\RSw{6.0cm}
\begin{figure}
\begin{center}
\includegraphics[width=\RSw,height=\RSh]{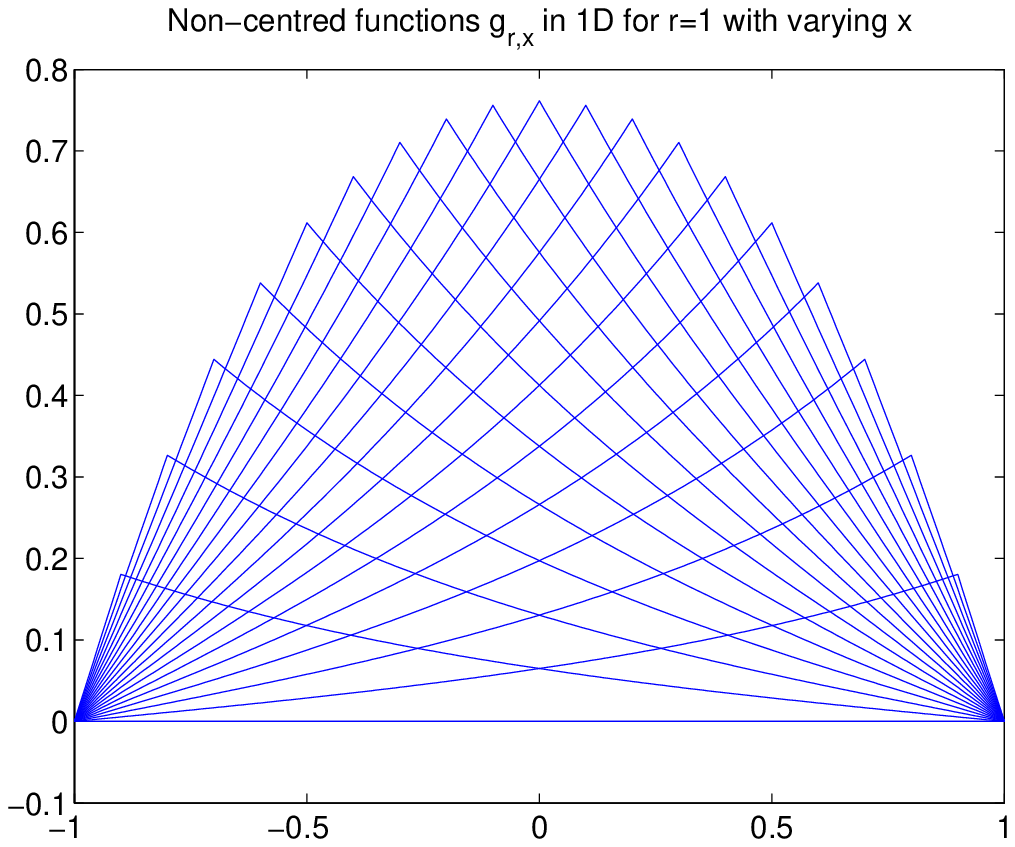} 
\caption{Kernel translates $g_{1,x}$ for the
  1D exponential kernel for varying $x$.
  \RSlabel{Fig1Dgrx}}
\end{center}  
\end{figure}
\begin{Theorem}\RSlabel{Thezeroinf}
If $\calk=\Wmd$, the kernel $K_r$ of \eref{eqKr}
is the reproducing kernel of the space $\Wmdr$.\qed
\end{Theorem} 
It is easy to generalize this to smoothly bounded
domains $\Omega\subset \R^d$
instead of the ball $B_r(0)$.
\biglf
There is a connection to Power Functions
on infinite point sets. We can define the Power Function for
infinite data outside the ball $B_r(0)$ as 
\bql{eqPowr}
P_{\|x\|\geq r}(x)=\sup\{v(x)\;:\;\|v\|_K\leq 1, v\in V_r\}.
\eq
Our standard criterion \eref{eqKvK0} for nonexistence of centralized bump
functions now has to be replaced by
$V_r=\{0\}$, and then the Power Function is zero.
\begin{Theorem}\RSlabel{ThePowBump}
In general,
$$
P_{\|x\|\geq r}(x)=\|g_{r,x}\|_K \fa x\in \R^d,
$$
and if there are bump functions,
the supremum in \eref{eqPowr}
is attained at $g_{r,x}(x)/\|g_{r,x}\|_K$. In particular,
$P_{\|x\|\geq r}(0)$ is attained at $b_r^*/\|b_r^*\|_K$, proving
\bql{eqUncRel}
P_{\|x\|\geq r}(0)=b_r^*(0)/\|b_r^*\|_K=\|b_r^*\|_K^{-1}.
\eq
\end{Theorem} 
\begin{proof}
  We assume that there are bump functions. Then  $g_{r,x}/\|g_{r,x}\|_K$
  is admissible, and
  $$
P_{\|x\|\geq r}(x)\geq g_{r,x}(0)/\|g_{r,x}\|_K=\|g_{r,x}\|_K.
$$
For all $v\in V_r$ we have $(v,g_{r,x})_K=v(x)$, and therefore
$$
|v(x)|\leq \|v\|_K \|g_{r,x}\|_K
$$
proves $P_{\|x\|\geq r}(x)\leq\|g_{r,x}\|_K$. In the nonexistence case,
all $g_{r,x}$ are zero, like the Power Function.
\end{proof}
The relation \eref{eqUncRel} was observed already in \RScite{schaback:2023-1}
as the extremal situation of a trade-off principle that relates small Power
Functions to large norms of bump functions.
\biglf
Theorem \RSref{ThePowBump} has an interpretation in Spatial Statistics.
If a random field $R$ on $\R^d$ has observations in all $x$ with
$\|x\|_2\geq r$,
$$
P_{\|x\|\geq r}(0)^2=\|b_r^*\|_K^{-2}
$$
is the variance of the Kriging predictor, i.e. the best linear unbiased
estimator from the data. For \M{} covariances generating $\Wmd$,
Section \RSref{SecBFSL} has shown that the variance behaves like $r^{2m-d}$.
Cases with nonexistence of bump functions are discouraged in Spatial
Statistics, because information on an infinite point set
like the complement of a ball implies total information. 
\biglf
Again, everything works the same for general domains $\Omega\subset\R^d$ with a
fixed interior point $x$. One has to project $K(x,\cdot)$
onto
$$
V_{\Omega}=\{v\in\calk\;:\;v(y)=0 \fa y \notin \Omega\}
$$
and renormalize the result to be 1 at $x$.
\section{Analytic Characterization}\RSlabel{SecAnalyChar}
Here, we focus on $\calk=\Wmd$ and want to
apply Real Analysis to find more specific
results on norm-minimal bump functions, including formulae ready for
computational implementation. Since the problem
of constructing norm-minimal bump functions can be written as a quadratic
optimization with infinitely many
constraints in an infinite-dimensional space, there is a variational
problem with Lagrange multipliers in the background, but we proceed directly to
the constructive solution.   
\biglf
As already stated,
the subspace $V_r$ of \eref{eqKvK0} is the closure of $C_0^\infty(B_r(0))$
under the native space norm. 
In case $\calk=\Wmd$, it is 
$\Wmdr$ 
in standard Sobolev space theory, and the boundary conditions are
well-known \RScite{schaback:2018-2}. By sources on trace theorems,
e.g. \cite[Thm, 10, Thm. 11]{suslina:2024-1}, the boundary
conditions for embedding $\Wmdr$ into $\Wmd$ consist of the classical
radial and normal derivatives 
$$
\gamma_{j,r}=\dfrac{\partial^j}{\partial \nu^j}|_{\partial B_r(0)},\;0\leq j<m-1/2
$$
whose extensions map $W_2^m(B_r(0))$ to $W_2^{m-j-1/2}(\partial B_r(0)),\;0\leq
j<m-d/2$. We integrate these over the boundary to define the functionals
$$
\lambda_{j,r}(f)=\displaystyle{\int_{\partial B_r(0)}\gamma_{j,r}(f)_{|_x}dx,
,\;0\leq j<m-1/2} 
$$
and the functions
$$
g_{j,r}(x):=\lambda^y_{j,r}K(x,y),
$$
where $\lambda_{j,r}$ acts on the variable $y$, as indicated by the
superscript. These are radial, i.e. rotationally invariant, because
the traces of $K(x,\cdot)$ on the boundary just rotate
with the direction of $x$, but the integral over the boundary stays the same.
\biglf
Next, we need the positive definite kernel matrix with entries
$$
\lambda^x_{k,r}\lambda^y_{j,r}K(x,y)=\lambda^x_{k,r}g_{j,r}(x),\;0\leq j,k<m-1/2
$$
and solve the $r$-dependent linear system
\bql{eqlKcllK}
g_{k,r}(0)=\lambda^y_{k,r}K(y,0)=\displaystyle{
  \sum_{j=0}^{j<m-1/2} c_j(r)\lambda^x_{k,r}\lambda^y_{j,r}K(x,y)  },
\;0\leq k< m-1/2
\eq
for functions $c_j(r),\;0\leq j<m-1/2$.
\begin{Theorem}\RSlabel{Thefinalrep}
  Norm-minimal compactly supported
  functions $b_r^*$ for Sobolev spaces $\Wmd$
  can be calculated via the radial function
\bql{eqgry}
g_r(y):=K(0,y)-\displaystyle{\sum_{j=0}^{j<m-1/2}c_j(r)g_{j,r}(y)}
\eq
in the above way, finally using \eref{eqbrgr}.
\end{Theorem}
\begin{proof}
By Theorem \RSref{TheChar} of Section \RSref{SecChar},
we need the projection $P_{\Wmdr}$ in $\Wmd$ onto $\Wmdr$ and the function
$$
g_r(\cdot)=P_{\Wmdr}K(0,\cdot)
$$
which is $b_r^*$ up to a factor by $b_r^*(x)=g_r(x)/g_r(0)$.
Due to \eref{eqvrgrvr}, the necessary and sufficient optimality conditions are
$g_r \in  \Wmdr$ and
$$
\begin{array}{rcl}
(K(0,\cdot) -g_r,v_r)_{\Wmd}&=&0 
\end{array} 
$$
$\fa r\in V_r=\Wmdr$. But \eref{eqgry} implies
$$
\begin{array}{rcl}
 (K(0,\cdot)-g_r,v_r)_{\Wmd}
  &=&
  \displaystyle{\left(\sum_{j=0}^{j<m-1/2}c_j(r)g_{j,r},v_r
    \right)_{\Wmd}}\\
  &=&
    \displaystyle{\sum_{j=0}^{j<m-1/2}c_j(r)(g_{j,r},v_r)_{\Wmd}}\\
  &=&
    \displaystyle{
      \sum_{j=0}^{j<m-1/2}c_j(r)(\lambda^y_{j,r}K(\cdot,y),v_r)_{\Wmd}}\\
  &=&
    \displaystyle{
      \sum_{j=0}^{j<m-1/2}c_j(r)\lambda^y_{j,r}v_r(y)}\\
  &=& 0  \fa v_r\in V_r.
\end{array}
$$
The system \eref{eqlKcllK} is the same as 
$$
\lambda^y_{k,r}g_r(y)=0,\;0\leq k<m-1/2,
$$
but since $g_r$ is radial, the boundary values
and radial derivatives are constant and therefore zero.
\end{proof}
This seems to be the first
case handling infinitely many data with infinitely many kernel
translates, in this case
placed on the $r$-sphere and treated in a rotationally symmetric way.
The definition of the functions
$g_{j,r}$ cares for orthogonality to $V_r$ and involves
all of these translates fairly.
The linear system \eref{eqlKcllK} has a different purpose: it cares
for the correct smoothness of the result on the boundary. If the $g_{j,r}$
and their derivatives are calculated on the boundary with sufficient accuracy,
the system can be set up and solved like any other Hermite interpolation.
\biglf
For domains $\Omega$ with a fixed interior point, the proof logic stays the
same, but the radiality arguments fail. Everything works as long as the trace
theorems are valid, but this fails for pathological subdomains.
\biglf
We add a remark on the background, connected to
the old theory of $L$-splines \RScite{schultz-varga:1967-1}.
The inner product
$(f,g)_K$ can be written as $(Lf, Lg)_{L_2(\R^d)}$
for a pseudodifferential operator
$L\;:\;\calk \to L_2(\R^d)$ defined by
$$
\widehat{Lf}(\omega)=\hat f(\omega)\sqrt{\hat K(\omega)}.
$$
Then the reproduction equation
$f(x)=(f,K(x,\cdot))_K=(Lf,LK(\cdot,x))_{L_2(\Omega)}$
is a way to define that $(L^*L)^yK(y,x)=\delta_x$
holds for all $x\in \R^d$, i.e. the kernel
is a {\em fundamental solution}. In the Sobolev case $\Wmd$, the operator
$L^*L$ is the classical elliptic differential operator $(Id-\Delta)^m$.
Compactly supported functions on a smooth domain $\Omega$ must then
obey the boundary  conditions for the Dirichlet problem for $L^*L$ 
on $\Omega$ with zero boundary conditions. Consequently, bump
functions may already be present in the literature
on elliptic PDE problems. Anyway, they are possibly useful
in the Method of Fundamental Solutions
\RScite{Bogomolny:1985-1,chen-et-al:2008-1}.
\section{Examples}\RSlabel{SecExa}
We first consider the simplest Sobolev space
$\Wmd$ for $m=1=d$ with the radial exponential
kernel $K(x,y)=\exp(-|x-y|)$. Then the 
functional
$$
\lambda_{0,r}(v)=v(-r)+v(r) \fa v\in W_2^1(\R^1)
$$
``integrates'' over the trace operator $T(v)=(v(r),v(-r))\in \R^2$.
With some explicit calculations omitted,
$$
\begin{array}{rcl}
g_{0,r}(x)&=&\displaystyle{\exp(-|r-x|) +\exp(-|-r-x|)}\\
g_r(x)&=&\exp(-|x|)-c_0\displaystyle{
  (\exp(-|r-x|) +\exp(-|-r-x|)}\\
c_0(r)&=&\dfrac{1}{\cosh(r)}\\
g_r(x)&=&\dfrac{ \sinh(r-|x|)}{\cosh(r)}\\
b_r^*(x)&=&\dfrac{g_r(x)}{g_r(0)}
=\dfrac{ \sinh(r-|x|)}{\sinh(r)}
\end{array} 
$$
where the third equation follows somewhat easier
from setting $g_r(x)=0$ for $x=\pm r$. This
is the situation of Figure \RSref{Fig1Dbumps}.
The corresponding general kernel translates
in the sense of \eref{eqKr} are in Figure \RSref{Fig1Dgrx},
obtained via projection of $K(x,\cdot)$ instead of $K(0,\cdot)$.
\def\RSh{6.0cm}
\def\RSw{6.0cm}
\begin{figure}
\begin{center}
\includegraphics[width=\RSw,height=\RSh]{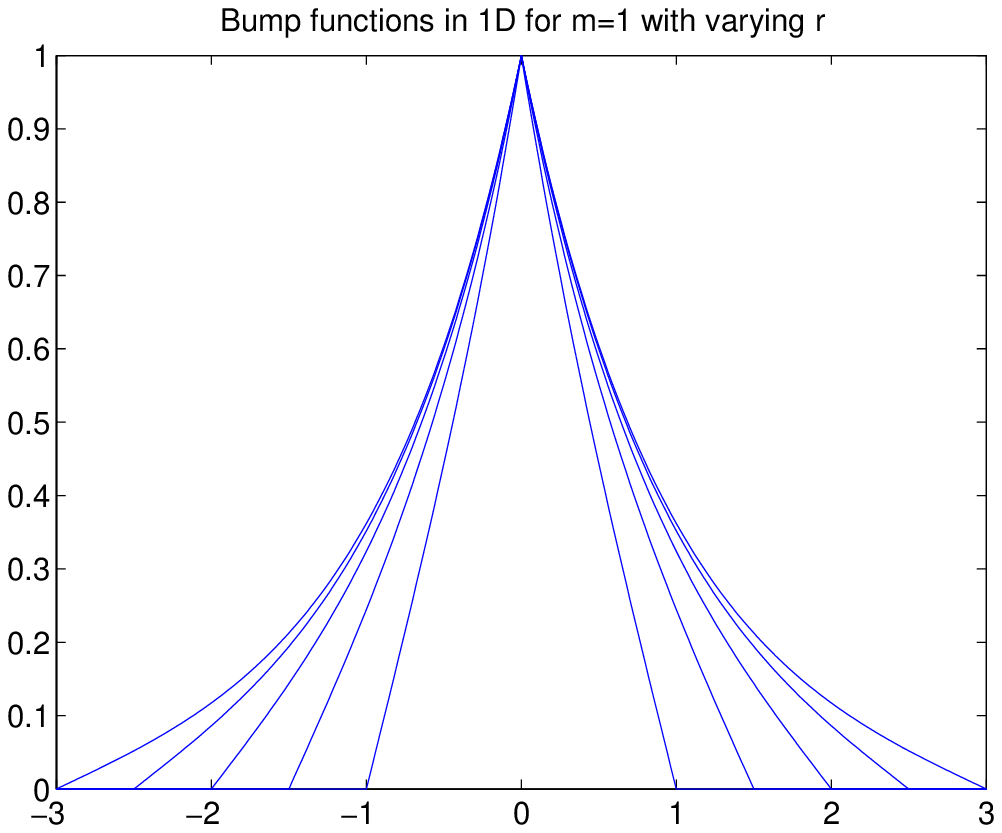} 
\includegraphics[width=\RSw,height=\RSh]{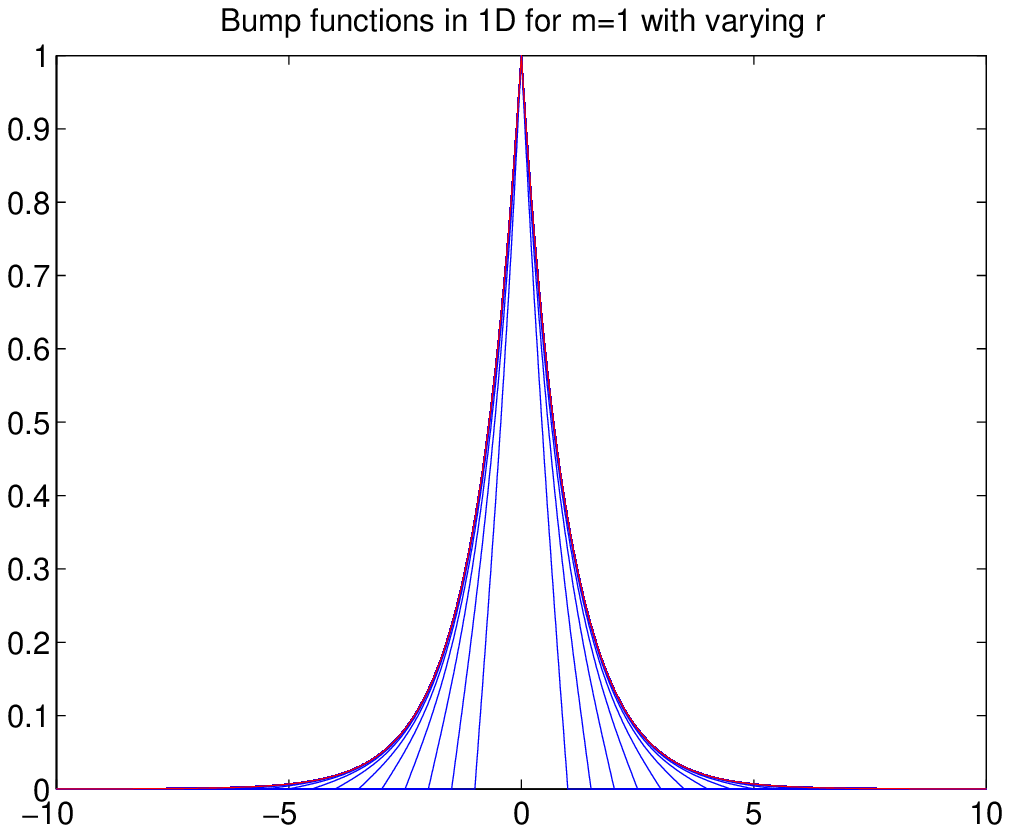}
\caption{Bump functions for the 1D exponential kernel for varying $r$.
  Left: $r=1,\,1.5,\,2,\,2.5,\,3$, right: $r=1,\ldots,50$. The red line
  is the kernel $K(0,y)=\exp(-|y|)$ occurring numerically in the limit.
  \RSlabel{Fig1Dbumps}}
\end{center}  
\end{figure}
\biglf
The same kernel works for all cases with $2m-d=1$, but things get much more
difficult for $d>1$, because we now have to integrate
over circles, using
$$
\lambda_{0,r}(f)=\displaystyle{\int_0^{2\pi}f(
  r\cos \varphi + r\sin(\varphi))d\varphi.   } 
$$
To keep things simpler, we shall use radiality whenever possible.
\biglf
We first consider $d=2,\, m=3/2$ and 
set
$x=r(\cos(\varphi),\sin(\varphi))$ and $z=(1,0)$
to get
\bql{eqxtz}
\|x-tz\|^2=(r\cos(\varphi)-t)^2+(r\sin(\varphi))^2
\eq
and
$$
\begin{array}{rcl}
g_{0,r}(tz)
&=&
\displaystyle{\int_0^{2\pi}   
  \exp\left(-\sqrt{(r\cos(\varphi)-t)^2+(r\sin(\varphi))^2}\right)d\varphi}\\
&=&
\displaystyle{\int_0^{2\pi}   
  \exp\left(-\sqrt{r^2+t^2-2r|t|\cos(\varphi)}\right)d\varphi}
\end{array} 
$$
because the cases $t$ and $-t$ have the same trace on the circle.
The ingredients for  \eref{eqlKcllK} are $g_{0,r}(0)=2\pi\exp(-r)$ and 
$$ 
\lambda^x_{0,r}\lambda^y_{0,r}K(x,y)=\lambda^x_{0,r}g_{0,r}(x)
=2\pi g_{0,r}(rz)
$$
because $g_{0,r}$ is constant on the boundary and equal to  $g_{0,r}(rz)$.
Therefore
$$
c_0(r)=\dfrac{\exp(-r)}{g_{0,r}(rz)}
$$
and we get
$$
g_r(x)=\exp(-\|x\|_2)-\dfrac{\exp(-r)g_{0,r}(x)}{g_{0,r}(rz)}.
$$
Figure \RSref{Fig1D2Dbumps} shows the norm-minimal bump functions (solid)
in the 1D (blue) and 2D case (red). The 2D case has a vanishing derivative
at zero, visible by zooming in, but still just continuity at $r$. But note that
functions in $W_2^{3/2}(\R^2)$ are only continuous, not continuously
differentiable. For comparison, the Wendland functions are added as dashed
lines. It is remarkable that in the 2D case the Wendland functions
have derivative discontinuities at zero,
while the bump functions have them at $r$. In case $r=10$, the exponential decay
of bump functions is apparent, while Wendland functions decay polynomially.
Smaller $r$ behave much like scaled versions of the case $r=1$.
\def\RSh{6.0cm}
\def\RSw{6.0cm}
\begin{figure}[hbt]
\begin{center}
\includegraphics[width=\RSw,height=\RSh]{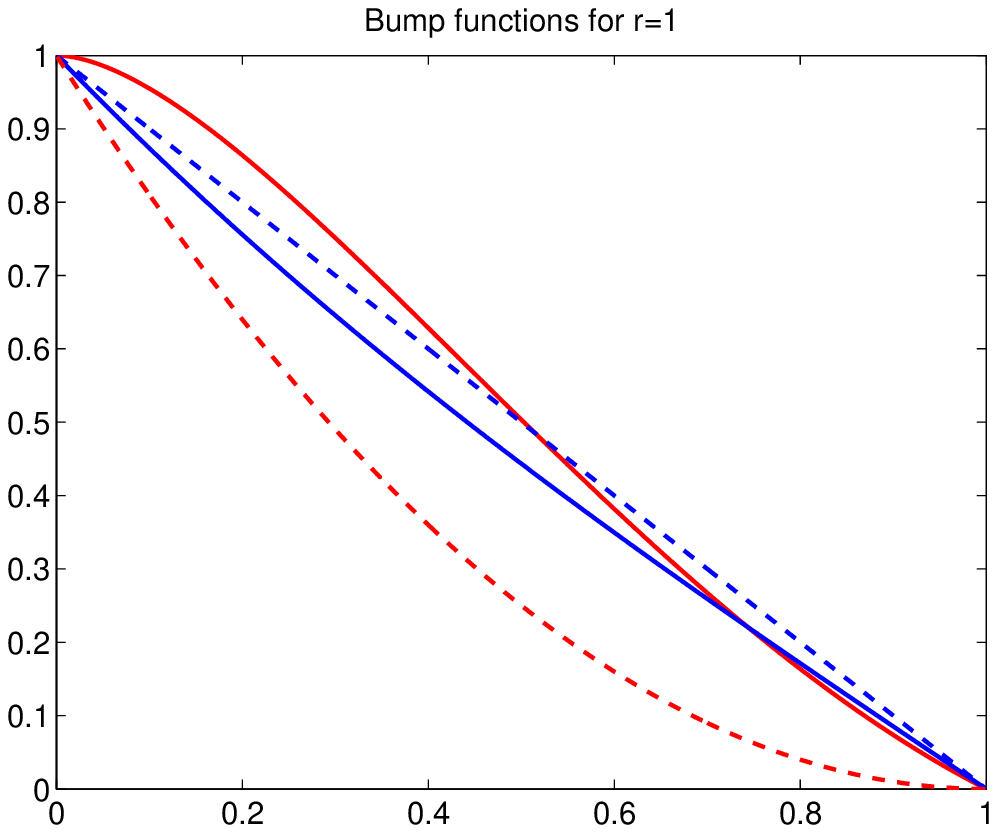} 
\includegraphics[width=\RSw,height=\RSh]{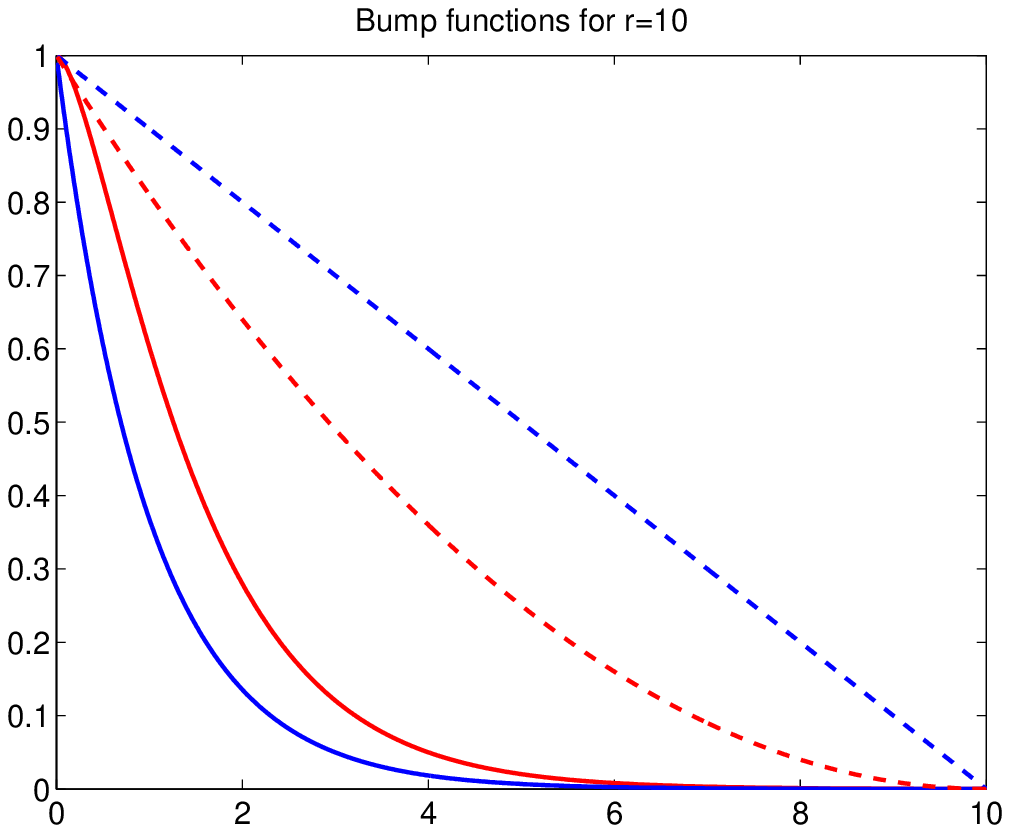}
\caption{Bump functions (solid)
  for the exponential kernel, $r=1$ (left) and
  $r=10$ (right). The red curves are for $W_2^{3/2}(\R^2)$, the blue curves for
  $W_2^1(\R^1)$. The corresponding Wendland functions are dashed. 
  \RSlabel{Fig1D2Dbumps}}
\end{center}  
\end{figure}
\biglf
Staying in the bivariate setting with larger $m$ just uses different kernels,
but now we get the linear system \eref{eqlKcllK} to solve. If we want
a derivative condition in $\R^2$, we can take $m=2.5$ or $m=2$.
When examining the case $m=2$, the function $g_{1,r}$ has a singularity
at the boundary, due to nonexistence of second derivatives in $W_2^2(\R^2)$.
For $m=2.5$, results are in Figure \RSref{Fig2Dbumps52}.
The corresponding Wendland function $\phi_{3,1}(t)=(1-t/r)_+^4(1+4t/r)$
is in $C^2$ and vanishes of third order at $r$. The norm-minimal
compactly supported function seems
to have the same smoothness at zero and at the boundary.
\def\RSh{6.0cm}
\def\RSw{6.0cm}
\begin{figure}
\begin{center}
\includegraphics[width=\RSw,height=\RSh]{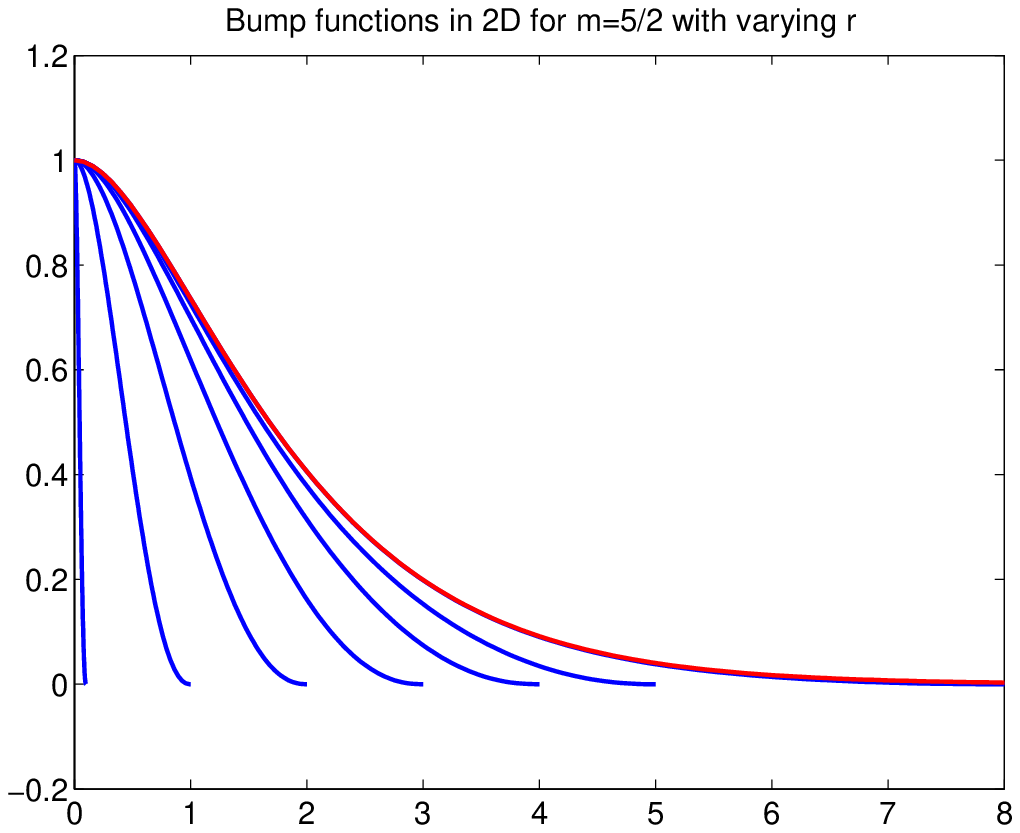} 
\includegraphics[width=\RSw,height=\RSh]{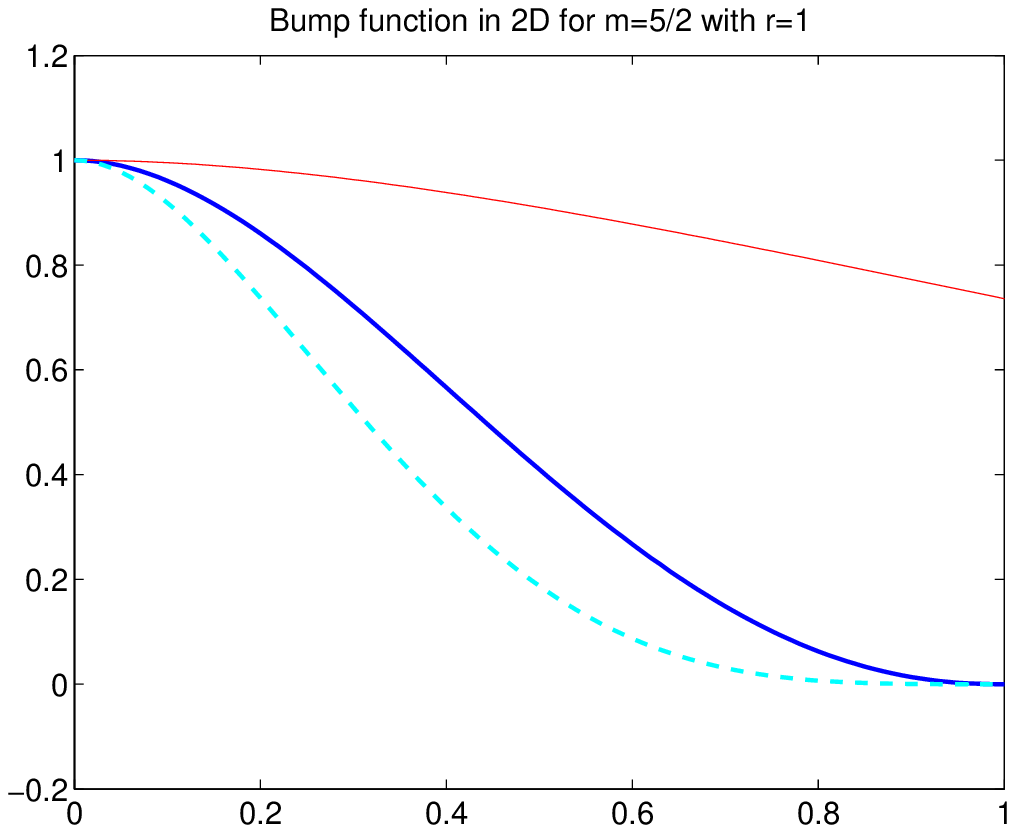}
\caption{Left: optimal bump functions
  for the 2D exponential kernel for $W_2^{5/2}(\R^2)$
  with varying $r$.
  The red line
  is the normalized kernel occurring in the limit. Right:
  the case $r=1$ with the Wendland function $\phi_{3,1}$ in cyan and dashed.
  \RSlabel{Fig2Dbumps52}}
\end{center}  
\end{figure}
\section{Summary, Conclusions, and Open Problems}\RSlabel{SecConc}
For a given support radius $r$, each Sobolev space $\Wmd$
with $m>d/2$ has a unique ``bump'' 
function with minimal norm, value 1 at the origin and vanishing
outside the interior of the ball $B_r(0)$ with radius $r$ around the origin.
Under all compactly supported functions with norm up to one and
value one  at the origin, it has smallest support.
And under all functions with support in $B_r(0)$ with norm up to one,
they attain the maximum value at zero, up to a factor.
\biglf
Such functions are an
intrinsic
and characteristic feature of the space and could be called $S$-splines.
\biglf
They are radial, and their radial univariate profile
can be precalculated to high accuracy by the construction
of Section \RSref{SecAnalyChar}. By the results on scaling
in Section \RSref{SecBFSL},
they can be downscaled to smaller support radii,
at asymptotically no loss. 
There are connections to upper bounds of convergence
rates of kernel-based interpolation, because they are Power Functions
for data outside a given ball. In Spatial Statistics using random fields
with \M{} kernels, they give the variance for Kriging estimation at zero
provided that there is full information outside the $r$-ball.  
\biglf
In Hilbert space terms, these functions 
are obtainable by renormalization
of the projection
of the reproducing kernel in $\Wmd$
onto $\Wmdr$. In terms of Real Analysis, they are computable by
applying the trace operator to the reproducing kernel under radial symmetry.
\biglf
For interpolation or approximation in Sobolev spaces,
their translates provide a compactly supported radial basis, leading to sparse
matrices, but there are no results yet on linear independence or
positive definiteness. In particular, their Fourier transforms are
not yet known, but they must have an optimality property as well.
The functions seem to be bell-shaped
\RScite{dyn-levin:1981-1} and pointwise decreasing for increasing $r$,
but this is still open. 
For use in meshless methods for PDE solving,
they are {\em shape functions} \RScite{belytschko-et-al:1996-1}
that deserve further investigation. Finally, they
may lead to new multivariate wavelet constructions, like many other
compactly supported functions.
\biglf
To make the numerical application of optimal compactly supported functions
easier, a follow-up paper should publish the
radial profiles in full computational accuracy. This
opens the way to various sparse
meshless methods for interpolation, approximation, and
PDE solving. Then it is interesting to see how much
bandwidth is needed to let them work 
at maximal possible convergence rate.
\biglf
The results of this paper should generalize to any kernel-based Hilbert space
with limited smoothness, and hopefully also to Beppo-Levi spaces
generated by conditionally positive definite kernels. Also,
the restriction to balls centred around the origin is easy to overcome,
losing radiality arguments. Then it is interesting to construct them on
tilings of the space, like finite elements. Their superposition will
stay in Sobolev space because of zero boundary conditions. 
\subsection*{Final Remarks}\RSlabel{SecFinRem}
There was no funding except the standard retirement program for professors
in the state of Lower Saxony, Germany.
\biglf
This work would not have been possible without the long-term friendship with 
Elisabeth Larsson (Uppsala) and Oleg Davydov (Gießen), leading to publications
\cite{larsson-schaback:2023-1,davydov-schaback:2019-1}. Bump functions were a
central tool in these papers.
\bibliographystyle{plain}

\end{document}